\documentclass{article}
\usepackage{amsbsy,amssymb,amscd,amsfonts,latexsym,amstext,delarray,
amsmath, diagrams}  \headsep=15pt
\usepackage[all,knot]{xy}
\def\piet{\pi_{1,\mbox{\'et}}}

\def\Un{ \mbox{Un} }

\def\bX{ \bar{X}}

\def\dim{{ \mbox{dim} }}
\def\Spec{{ \mbox{Spec} }}

\def\ra{{ \rightarrow }}
\def\da{{ \downarrow }}

\def\a{{ \alpha }}
\def\b{{ \beta} }

\def\Un{ \mbox{Un} }
\def\Vect{ \mbox{Vect} }
\def\hra{{ \hookrightarrow }}
\def\da{{ \downarrow }}

\def\C{{ \mathbb{C} }}
\def\R{{ \mathbb{R}}}
\def\bs{ \backslash}

\def\G{{ \Gamma }}
\def\Gal{{ \mbox{Gal} }}

\def\bQ{\bar{\Q}}

\def\Z{{ \mathbb{Z}}}

\def\bq{\begin{quote}}
\def\eq{\end{quote}}
\def\Aut{ \mbox{Aut}}
\def\Isom{ \mbox{Isom} }

\newtheorem{thm}{Theorem}[section]

\def\Q{\mathbb{Q}}

\def\invlim{\varprojlim}

\def\be{\begin{equation}}
\def\ee{\end{equation}}

\def\tX{\tilde{X}}

\def\tV{\tilde{V}}

\def\bE{ \bar{E}}

\def\Pet{\pi_1^{et}(\bX;b,x)}
\def\piet{\pi_1^{et}(\bX,b)}
\def\b{ \bar}

\def\hT{\hat{T}}

\def\XC{ X(\C)}

\def\tbE{\tilde{\bar{E}}}

\def\a{\alpha}
\def\Ph{P^{\wedge}}
\def\Gh{G^{\wedge}}
\def\bd{\begin{diagram} }
\def\ed{ \end{diagram}}
\def\tXC{ \widetilde{\XC}}
\def\tbX{ \widetilde{\bar{X}}}
\def\b{\beta}
\def\tb{\tilde{b}}

\def\bs{\begin{slide}}
\def\es{\end{slide} }

\def\Cov{\mbox{Cov}}
\def\bd{ \begin{diagram}}
\def\ed{\end{diagram} }
\def\tbX{\tilde{\bar{X}}}
\def\Cov{\mbox{Cov}}
\def\ux{\underline{x}}
\def\Loc{\mbox{Loc}}

\def\Pux{\pi_1^{u,\Q_p}(\bX;b,x)}
\def\piu{\pi_1^{u,\Q_p}(\bX,b)}
\def\pia{\pi_1^{alg,\Q_p}(\bX,b)}
\def\Px{\pi_1(\XC;b,x)}
\title{
Fundamental groups and Diophantine geometry}
\author{Minhyong Kim}
\begin{document}
\maketitle
{\em Colloquium lecture, Leeds, January 2008}
\medskip

To work our way towards the very canonical but rather difficult relationship between the notions
appearing in the title, it is appropriate to review briefly the
classical problems that make up the background of our study, and whose
importance will be initially regarded as self-evident.
Thus, we are given a polynomial $$f(x_1,x_2,\ldots, x_n)$$ whose coefficients
will be assumed to be in $\Z$ for the sake of simplicity.
The set of solutions to the equation
$$f(\ux)=0$$
can be considered in any number of different environments
such as
$$\Z,\  \Z[1/62], \ \Q,\  \Z[i],\  \Q[i], \ldots,\  \Q[i,\pi], \ldots, \ \R,\  \C,\  \Q_p, \ \C_p, \ldots$$
In recent decades, the designation of the  equation as Diophantine has not been a reference to any particular
property of the equation itself, but rather calls attention to our primary focus
on contexts closer to the beginning of the list, although how far we might extend the scope
is better left undetermined.  In any case, there
are famous results corresponding to different lines of demarcation, such as the one that says
$$x^n+y^n=z^n$$
has only the obvious solutions in $\Z$ as long as $n\geq 3$, or  where
$$f(x,y)=0$$
for a generic $f$ of degree at least 4 has only finitely many solutions in
$\Q(i,\pi,e)$.

Elementary {\em coordinate geometry} can be brought to bear on some such questions as a potent tool for describing solution sets, or
least for generating solutions.
A simple but already interesting case is a quadratic equation in two variables,
say
$$x^2+y^2=1.$$
By visualizing the real solution set as a circle, we might come upon the
idea of considering the intersections with lines that pass through
the specific point $(-1,0)$, where the set-up has already encouraged us casually to
refer to a solution  using  geometric terminology. The lines are described
 using equations $y=m(x+1)$ for various $m$ whereby  algebraic substitution leads to the constraint
$$x^2+(m(x+1))^2=1$$
or
$$(1+m^2)x^2+2m^2x+m^2-1=0.$$
A deeper connection to algebra comes from
the observation that one solution $x=-1$ is already rational, so that whenever the slope
$m$ is rational, the other solution is also
bound to be rational. As we vary $m$, we can generate thereby {\em all } the other rational
solutions to the equation, for example,
$(-99/101, 20/101)$ corresponding to $m=10$.
It seems that the visually compelling nature of the solution set
in a sufficiently big field provides valuable insight into finding solutions in
much smaller fields.
Incidentally, I'm sure you're aware also that this procedure
leads to the famous {\em Pythagorean triples} involved in equations like
$$99^2+20^2=101^2.$$

The elementary elegance of the method described becomes progressively harder to retain with the
increasing complexity of the problem, measured, for example, by the
degree of the equation. Nevertheless, it is instructive to
consider one  example of degree 3:
$$x^3+y^3=1729.$$
One verifies with the help of Ramanujan that
$(9,10)$ is a solution, so the case of
the circle might motivate us to consider lines through it. Unfortunately, the previous argument for
the rationality of intersection points fails as the associated constraint
becomes cubic. But if we want to start out generating just {\em one}  other solution, a  more subtle
idea is to consider  the tangent line to the real curve
at the point $(9,10)$, because then, the corresponding cubic equation will have
9 as a {\em double} root. To spell this out, calculate the equation of
the tangent line,
$$81(x-9)+100(y-10)=0$$
or
$$y=(-81/100)x+1729/100,$$
and substitute to obtain the equation
$$x^3+((-81/100)x+1729/100)^3=1729.$$
We have arranged for $x=9$ to be a double root, and hence,  the remaining root is forced
to be rational. Even by hand, you can (tediously) work  out
the resulting rational point to be
$$(-42465969/468559, 24580/271).$$
Repeating the procedure with the points that are successively obtained thus
actually provides us with infinitely many rational solutions.
Here, you must pause to consider the possibility that  repetition
will just move us (quasi-)periodically around finitely many points,
but there is a well-known theorem of Nagell and Lutz that tells us this cannot
happen given the denominator of the solution at hand.

Geometric techniques of the same general flavor can be made considerably
more sophisticated, with nice applications to varieties of simple type as might be
defined by equations of low degree in a greater number of variables.
But in the present lecture we wish to  explain the important
conceptual shift that occurred in the 1960's, whereby Diophantine problems acquired
an {\em intrinsically } geometric nature by way of two foundational
ideas of Grothendieck.

The first one, elementary in comparison to the second, associates to
the polynomial
$f(\ux)$ the ring
$$A:=\Z[\ux]/(f(\ux)).$$
This leads to a natural correspondence between
solutions
$(b_1, \ldots, b_n)$ of
$f(\ux)=0$ in a ring $B$, and ring homomorphisms
$$A\ra B$$
That is, an {\em arbitrary} n-tuple
$\underline{b}=(b_1,\ldots, b_n)$ determines a ring homomorphism
$\Z[\ux]\ra B$ that sends $x_i$ to $b_i$, which factors through the
quotient ring $A$ exactly when $\underline{b}$ is a zero of $f(\ux)$.
The spatial intuition is supposed to arise from the idea that a
commutative ring $R$ with 1 can be viewed as the ring of functions on
a space, its {\em spectrum} $$\Spec(R),$$
whose underlying set consists of the prime ideals of $R$. This correspondence reverses arrows reflecting
the intuition that
a map of spaces should pull functions backwards by composition. Thus, the
solutions in $B$ of
$f(\ux)=0$ come into bijection with the set of
maps
$$\Spec(B) \ra X:=\Spec(A),$$
conventionally denoted by $$X(B).$$
Even before considering such difficult maps, it is pleasant to note that
an obvious map
$$\begin{array}{c}X \\
\da \\
\Spec(\Z)\end{array}$$
corresponds to the inclusion
$$\Z \ra A=\Z[\ux]/(f(\ux))$$ using which we think of
$X$ as a fibration over $\Spec(\Z)$.
Then the solutions in $\Z$, the elements
of $X(\Z)$, are precisely the {\em sections}
$$
\xy
(0,22)*{X};
(0,-2)*{\Spec(\Z)};
(0,20)*{}="A";
(0,0)*{}="B";
(3,0)*{}="C";
(3,20)*{}="D";
{\ar@{->}"A";"B"};
{\ar^{P}@/_1pc/"C";"D"};
\endxy
$$
of the fibration.
The remarkable upshot of this formulation is that the
study of solutions to equations is subsumed into the study of
maps whose very nature compels us to consider as the most basic  in all of  mathematics.
This perspective is of late provenance in the theory of Diophantine equations, but still
 provides at this point its most fundamental justification.

The second idea involves a sophisticated construction whereby
spaces like $\Spec(\Q)$ or $\Spec (\Z)$ are endowed with very non-trivial
topologies that go beyond  scheme theory (by which we mean the global theory of such spectra).
 We will not review the precise definitions in this summary, since it appears by now
 well-known that a Grothendieck topology on an object $T$ allows open sets to be
certain maps with range $T$ from domains that are not necessarily  subsets of $T$.
On a `usual' topological space, one could make the  topology finer by
allowing as open sets maps $$U\ra T$$ that factorize as
$$U\hra V \ra W \hra T$$
where $W\hra T$ is an open subset, $V\ra W$ is a covering space,
and $U\hra V$ is an open embedding. An open covering then
is a collection $\{U_i\ra T\}_{i\in I}$ of such maps
with the property that the union of the images is $T$.
But this does not give anything essentially new. By definition each such $U\ra T$ is a local homeomorphism,
so that coverings by families of usual open subsets is co-final
among all such exotic open coverings. That is to say, any covering
$\{U_i\ra T\}_{i\in I}$ in the generalized sense has a refinement
$$\{V_{ij} \hra T\}$$ where each $V_{ij}\hra T$ is an open embedding
that factors through one of the $U_i$:
$$V_{ij} \ra U_i \ra T.$$
This fact induces an equivalence of categories between the category of usual sheaves
and sheaves in this refined topology.

However, in algebraic geometry, there are many maps that behave formally
like local homeomorphisms without actually being so. These are the so-called
{\em \'etale maps} between schemes. A nice and fairly general class of examples arise from
maps
$$\Spec(B) \ra \Spec(A)$$ corresponding to  maps of rings
$A \ra B$ where $B$ has the form $$A[x]/(f(x))$$for
a monic polynomial $f(x)$. The constraint we wish to impose is that the fibers of
$\Spec(B)$ over $\Spec(A)$, which have the form
$$\Spec(k[x]/(\bar{f}(x)))$$
for residue fields $k$ of $A$,
 should have the same number of elements, indicating an absence
of ramification.
For this, we need to prevent $f(x)$ from having multiple roots
in any such residue field. This amounts to the condition that
$f(x)$ and $f'(x)$ should not have common roots point-wise,
or that the discriminant of $f$ should be a unit in $A$.
The obvious map
$$\Spec(\C[t][x]/(x^2-t))\ra \Spec(\C[t]),$$
is not \'etale,
the discriminant of $x^2-t$ being the non-unit $4t$,
while
$$\Spec(\C[t,t^{-1}][x]/(x^2-t))\ra \Spec(\C[t,t^{-1}]),$$
is \'etale.

Allowing \'etale maps as open subsets gives a genuinely richer topology
to a scheme than the Zariski topology.
The connected \'etale coverings of $\Spec(\Q)$, for example,
are maps
$$\Spec(F) \ra \Spec(\Q),$$ where $F$ is a finite field extension
of $\Q$. For $\Spec(\Z)$, one can construct an open covering using
the two maps
$$\Spec(\Z[i][1/2]) \ra \Spec(\Z)$$
and
$$\Spec(\Z[(1+\sqrt{-7})/2][1/7]) \ra \Spec(\Z).$$

The (co-)homology theory associated to sheaves in the
\'etale topology has been fabulously applied to the
arithmetic geometry of schemes in the past many decades, with results well-enough known
not to require a separate survey.
Less known perhaps, is that Grothendieck's exotic topologies can
also lead to interesting {\em homotopy} groups, whose structures
are only recently being probed at any depth.
One such direction is
the {\em motivic homotopy theory} of Voevodsky,  about which we
will say nothing. The emphasis here instead is on rather recent
developments in a somewhat older homotopy theory belonging to the {\em \'etale fundamental group}
and its variations.
In particular, we will focus exclusively on the
application of the theory to Diophantine problems.

 The beginning point is surprisingly elementary, wherefrom
  the theory obtains a substantial portion of its charm.
  Let therefore $X$ be a variety defined over $\Q$ and
 $G=\pi_1(\XC, b)$ the usual topological fundamental group
of the space obtained from the complex points of $X$.
For any point $x\in \XC$, we can also consider the
homotopy classes of paths
$$\pi_1(\XC;b,x)$$ from $b$ to $x$.
Then $\Px$ has the natural structure of a principal $G$-bundle, or a {\em $G$-torsor}, in that
$G$ naturally acts on $\Px$ via composition of paths, and the choice of any
$p\in \Px$  induces a bijection
$$G\simeq \Px$$
$$g\mapsto pg$$
via the action. Since this principal bundle lives on a topological
point, of course it is trivial. However, we see even here that the
{\em variation} of $\Px$ is $x$ is not at all trivial in general.
That is to say, the triviality of the individual $P_x$
is not different from the triviality of the fibers of even a complicated
vector bundle.
To be more precise on this point, choose a pointed universal covering space
$$f:(\widetilde{\XC}, \tb) \ra (\XC,b).$$
Then lifting of paths
determines  natural bijections
$$\Px\simeq \tXC_x$$
between homotopy classes of paths and the fibers of the universal covering
space. In fact, it is natural to {\em construct}
$\tXC$ as
$$\cup_x \pi_1(\XC;b,x)$$
topologized so that
the obvious projection that takes $\pi_1(\XC;b,x)$ to $x$ is a local homeomorphism.
In any case, we see thereby that the principal bundles in question form
the fibers of a map $$f:(\widetilde{\XC}, \tb) \ra (\XC,b)$$
that can be highly non-trivial. In fact, we will see that
the lack of a {\em canonical} isomorphism $G\simeq \Px$  is the essential
ingredient underlying our ability to
endow $\Px$ with a genuinely non-trivial structure of a principal
$G$-bundle within suitably enriched contexts.

 As far as Diophantine problems are concerned, we will of course be interested in the
 situation where $b$ and $x$ are both rational points in $X(\Q)$.
 As it stands, the principal bundles $\Px$ cannot
 pick out such special points as being
 different in any way from generic points.
 There are several ways to remedy this, of which
 the (ostensibly) easiest one to explain is the passage to the pro-finite completion.
 That is, define
 $$\Gh:=\invlim_{N\lhd G, [G:N]< \infty} G/N$$
 and
 $$\Ph:=\invlim_{N\lhd G, [G:N]< \infty} P/N$$
 for any principal $G$-bundle $P$.
 Then the basic and remarkable fact is that
 $\Gh$ is a sheaf of groups on the \'etale topology of $\Spec(\Q)$
 while $\Px^{\wedge}$ is a principal bundle for $\Gh$ in this topology.
 This statement is demystified just a little bit
 by recalling that a sheaf on $\Spec(\Q)$ is simply a set equipped with
 a continuous action of $\G=\Gal(\bQ/\Q)$. Nevertheless, it remains to
 see that the Galois group will indeed act on an object that arose thus out of
 ordinary topology.

Accounting for the action is an isomorphism
 $$\pi_1(\XC,b)^{\wedge}\simeq \pi^{et}_1(\bX,b)$$
 where $$\bX=X\times_{\Spec(\Q)}\Spec(\bQ)$$ is $X$ regarded as a variety
 over $\bQ$, while $\pi_1^{et}$ refers to the {\em pro-finite \'etale fundamental group}.
 It is the latter object on which $\G$ will act naturally.

 The definition will be reviewed after a brief return to usual topology.
For a manifold $M$ and an element $b\in M$, the fundamental group
$\pi_1(M,b)$ of $M$ with base-point $b$ can be defined in at least two different
ways avoiding  direct reference to topological loops.
One way is to
note first that a loop $l$ acts naturally on the fiber over $b$ of any covering space
$N\ra M$ of $M$ using the monodromy of a lifting $\tilde{l}$ of $l$ to $N$:
$$l_N:N_b \simeq N_b$$
This bijection is compatible with composition of loops on the one hand,
and with maps between covering spaces, on the other. That is,
$(ll')_N=l_N\circ l'_N$, and if
$f:N \ra P$ is a map of covering spaces, then
$$f\circ l_N=l_P\circ f$$ as maps from $N_b $ to $P_b$.
It is something of a surprise that the only way to give such
a compatible collection of automorphisms is in fact using an element
of the fundamental group. The concise way to state this
is via the
functor
$$F_b:\Cov(M) \ra \mbox{Sets}$$
that associates to each covering $N$ its fiber $N_b$ over $b$.
Then the fact in question is that
$$\pi_1(M,b)\simeq \Aut(F_b)$$
with the $\Aut$  understood in the sense of
invertible natural transformations of a functor.

Now given a variety $V$, we can use this approach to {\em define}
the \'etale fundamental group simply by changing the category of coverings. So
we let $$\Cov^{et}(V)$$ be the finite \'etale covers of
$V$ and, for any point $b\in V$, consider the functor
$F^{et}_b$ that takes $W \ra V$ to the fiber $W_b$.
Then
$$\pi_1^{et}(V,b):=\Aut(F^{et}_b)$$
Similarly,
$$\pi_1^{et}(V;b,x):=\Isom(F^{et}_b,F^{et}_x)$$
These superb definitions have been around at least since the 1960's, but it is rather striking
that variation of the base-point has not been really attended to until
fairly recently. The primary impetus for a serious reassessment  appears to have
come from the interaction with the Hodge theory of the fundamental group.

Nevertheless, constructions of the same general nature have now become commonplace in mathematics,
 the best known  being associated to the notion of
a Tannakian category, whereby the automorphisms of suitable functors
defined on agreeable categories give rise to group schemes.
Here we will content ourselves with mentioning two more examples.
Fix a non-archimedean completion $\Q_p$ of $\Q$ and
consider the category
$$\Loc^{et}(V,\Q_p)$$ of locally constant sheaves of finite-dimensional $\Q_p$-vector
spaces on $V$ considered in the \'etale topology.
There is still a fiber functor
$$F^{alg}_b:\Loc^{et}(V,\Q_p) \ra \Vect_{\Q_p},$$
now taking values in $\Q_p$-vector spaces, that associates to each sheaf its stalk at $b$.
(In comparing with the previous situation, it would be useful for the audience to have some intuition for the notion that a locally constant
$\Q_p$-sheaf is a `linearized' version of a covering space.)
Now define
$$\pi_1^{alg,\Q_p}(V,b):=\Aut^{\otimes}(F^{alg}_b),$$ the {\em $\Q_p$-pro-algebraic completion}
of $\pi_1^{et}(V,b)$. The $\otimes$ in the superscript refers to the fact
 that the automorphisms are required to be compatible not just with the
 morphisms in the category, but also the tensor product structure.
 As the name suggests, it is a pro-algebraic group over $\Q_p$.

When we replace all local systems by unipotent ones, i.e., those that
admit
a filtration
$$L=L^0\supset L^1 \supset \cdots L^n\supset L^{n+1}=0$$
such that each quotient
$L^i/L^{i+1}$ is isomorphic to a direct sum of the constant sheaf
$\Q_p$, one again gets a category $\Un^{et}(V, \Q_p)$ of the right sort to which one
can restrict the previous fiber functor
$$F^u_b:\Un^{et}(V,\Q_p) \ra \Vect_{\Q_p}.$$
The {\em $\Q_p$-pro-unipotent completion} of the \'etale fundamental group is then defined
as
$$\pi_1^{u,\Q_p}(V,b):=\Aut^{\otimes}(F^u_b)$$
In both settings, there are still torsors of paths
$$\pi_1^{alg,\Q_p}(V;b,x):=\Isom (F^{alg}_b,F^{alg}_x)$$
and
$$\pi_1^{u,\Q_p}(V;b,x):=\Isom (F^u_b,F^u_x)$$

It is natural to regard such definitions with a degree of suspicion, since
 not having loops to visualize may make them
 seem entirely intractable. The situation is somewhat ameliorated
through the intermediary of a universal object, which we describe in detail only
for the full pro-finite \'etale fundamental group.
Because $\Cov^{et}$ consists of {\em finite} covering
spaces, it may not be possible to find a single universal object inside the category. However, it is possible to construct
a pro-object that performs the same role. This is a compatible system
$$\tV=\{V_i\}_{i\in I}$$ of finite \'etale coverings
$$V_i\ra V$$ indexed by some filtered set $I$, having the following universal
property: If we choose $\tb=(b_i) \in \tV_b$, the pair
$(\tV,\tb)$ is universal among pointed pro-covering spaces, in that
any finite \'etale pointed covering $(Y,y) \ra (V,b)$ fits into
a unique commutative diagram

$$\begin{diagram}
(\tV,\tb)& \rTo & (Y,y) \\
 & \rdTo & \dTo\\
 & & (V,b) \end{diagram}$$
This means that there is some index
$i$ and a commutative diagram
$$\begin{diagram}
(V_i,b_i)& \rTo & (Y,y) \\
 & \rdTo & \dTo\\
 & & (V,b) \end{diagram}$$
In this situation, once again we have essentially tautological isomorphisms
 $$\pi^{et}_1(V,b)\simeq \tV_b$$
 and
 $$\pi^{et}_1(V;b,x)\simeq\tV_x,$$
 where the fibers are also projective systems of points.

 When $V=\bX$ for a variety $X$ defined over $\Q$ and the base-point
 $b$ is in $X(\Q)$,
 then the entire pro-system $$\tbX \ra \bX$$ comes from a system
 $$\tX\ra X$$ defined over $\Q$ and we can choose
 $\tb\in \tbX$ as well to come from  a rational point
 $\tb \in \tX$. The isomorphisms
  $$\pi^{et}_1(\tbX;b,x)\simeq\tbX_x$$
  then are compatible with the action of $\G$.
  The sheaves on $\Spec(\Q)$ obtained thereby have also
  a harmonious description in terms of the map
  corresponding to a rational point.
  The point is that the map
  $$\tX\ra X$$
  is a pro-sheaf of sets in the \'etale topology of $X$.
  Then given any point
  $$x:\Spec(\Q) \ra X,$$
  we get the sheaf
  $$x^*(\tX)$$
  on $\Spec(\Q)$, which is nothing but
  $\pi^{et}_1(\tbX;b,x)$.

  We illustrate this construction with the example of
$(\bE,0)$, an elliptic curve with origin over $\Q$.
 Let $$E_n \ra E$$ be the covering space given
 by $E$ itself with the multiplication map
 $$[n]:E \ra E.$$
 Then the system
 $$\begin{diagram}
 (\tilde{\bE}, \tilde{0}):=\{(\bE_n,0)\}_n & \rTo &(\bE,0)
 \end{diagram}$$
 is a universal pointed covering space.
Thus,
 for $(\bE,0)$,
 $$\pi^{et}_1(\bE,0)\simeq \hT(E)$$
 and an element of the fundamental group is just
 a compatible collection of torsion points of $E$.
That is to say, the Galois action on $\pi^{et}_1(\bE,0)$
is the well-known action on the Tate module of $E$.
 Similarly,
 $$\pi^{et}_1(\bE;0,x)\simeq \tbE_x$$
 consists of compatible systems of division points of
 $x$.

 A notable fact that emerges from this description is that if we take into
 account the Galois action, it is
 no longer possible to trivialize the
 torsor in general, even point-wise.
That is, there will often be
 no isomorphism between
 $\pi^{et}_1(\bX,b)$ and $\pi^{et}_1(\bX;b,x),$
 reflecting the fact that the \'etale topology has a very rich structure
 even on a point.
 In the case of
 $(E,0)$, if there were an isomorphism
 $$\pi^{et}_1(\bE,0) \simeq \pi^{et}_1(\bE;0,x)$$
 then there would be a Galois invariant element of
 $$\pi^{et}_1(\bE;0,x)\simeq \tilde{\bE}_x.$$
 In particular,
 for any $n$, there would be a rational point $x_n$ such
 that $nx_n=x$, which is not possible for $x\neq 0$ by a theorem of Mordell.

 To summarize, given a variety $X/\Q$ with a fixed rational point
 $b\in X(\Q)$, we are associating to each other point $x\in X(\Q)$
 a principal bundle $\Pet$ for $\piet$
 on the \'etale topology of $\Spec(\Q)$.
 This information can be organized using  a standard classifying space of sorts for
 principal bundles. That is, given a principal bundle $T$, one can choose
 a point $$t\in T$$ and examine the action of $\G$ on that point.
 For each $g \in \G$, $g(t)$ will be related to $t$ by an element
 $l_g\in \piet$, that is, $$g(t)=tl_g.$$
 The map $$g\mapsto l_g$$ obtained thereby is a 1-cocycle
 $$c_t:\G \ra \piet,$$  that is, a continuous map that satisfies
 $$c_t(g_1g_2)=c(g_1)g_1(c(g_2)).$$
 If we denote the set of such cocycles by $$Z^1(\G, \piet),$$ then
 $\piet$ acts on it according to
 $$lc(g):=g(l^{-1})c(g)l$$ and a different choice of $s\in T$ will lead to
 a cocycle $c_s$ lying in the same orbit as $c_t$. Denote by
 $$H^1(\G,\piet):=\piet\backslash Z^1(\G, \piet)$$
 the orbit set, so that the torsor $T$ determines a
 class $$[T]=[c_t]\in H^1(\G,\piet).$$ This cohomology set in fact classifies
 such torsors so that we have defined a map
 $$\bd X(\Q) &\rTo & H^1(\G,\piet)\ed$$
 $$\bd x &\rMapsto & [\Pet]\ed$$
 to a classifying space that can be thought of as an {\em \'etale
 period map}.
 In his famous letter to Faltings, Grothendieck formulated  the
 hope of studying Diophantine problems using this map.
  (He did not express matters using torsors, but rather, splittings of
 a certain canonical sequence of fundamental groups, in order to better harmonize the
 discussion with his general program of {\em anabelian geometry}.)

Unfortunately, it seems at present that the set  $H^1(\G,\piet)$ has too little
structure to study in a comprehensible manner. It should be obvious, meanwhile,
that an entirely analogous construction can be carried out with
$\pi_1^{alg,\Q_p}(\bX,b)$ or with
$\pi_1^{u,\Q_p}(\bX,b)$. For reasons that are somewhat technical to discuss in a short survey,
$\pi_1^{alg,\Q_p}(\bX,b)$ does not afford much advantage at present over
$\piet$. The unipotent completion, on the other hand, has been exploited to a certain extent
in the study of Diophantine sets. The key difference from the
other cases has to do with the relative ease of accessing information about
$$H^1(\G, \piu),$$
or rather, a slight improvement of this set.
Let $S$ be the set of primes of bad reduction for $X$, and denote by
$X(\Z_S)$ the set of points in the ring $\Z_S$ of $S$-integers, where the integrality
is defined in terms of a suitably good model. (Note that if $X$ is compact, then
the integral points are the same as rational points.)
The first point of note is that the map
$$X(\Q) \ra H^1(\G, \piu)$$
$$x\mapsto [\Pux],$$
when restricted to the integral points,
factors through a natural subspace
$$\begin{diagram}
X(\Z_S)&\rTo &H^1_f(\G, \piu)\subset H^1(\G, \piu)
\end{diagram}$$
corresponding to  local conditions satisfied by the
torsors $[\Pux]$, such as being unramified away from the primes of bad reduction
and $p$, and having a `crystalline' nature at $p$. This last  condition arises
from the {\em $p$-adic Hodge theory} of the non-archimedean variety $$X\times_{\Spec(\Q)}\Spec(\Q_p)$$
 that exerts a useful influence on $\piu$.
 In fact, these  conditions are meaningless for
$H^1(\G, \piet)$ and  quite difficult to analyze for $H^1(\G, \pia)$.
The advantage of considering them  in the unipotent setting is that
the subspace
$H^1_f(\G, \piu)$ becomes canonically equipped with the structure of a pro-algebraic variety.
In fact, for various quotients $[\piu]_n$ of $\piu$ modulo its descending central series,
 the sets
$$H^1_f(\G, [\piu]_n)$$ have natural structures of algebraic varieties over
$\Q_p$ that fit into a tower:
$$
\begin{diagram}[height=1.5em]
 &  & \vdots \\
  & \vdots &  &H^1_f(\G,[\piu]_4) \\
  &\ruTo(2,6)&  &\dTo \\
 & & &H^1_f(\G,[\piu]_3)\\
 &\ruTo(2,4) & &\dTo\\
 &  & &H^1_f(\G,[\piu]_2) \\
 & \ruTo  & &\dTo\\
X(\Z_S) & \rTo& &H^1_f(\G,[\piu]_1) \\
\end{diagram}$$
refining the map at the bottom (which has a classical interpretation in Kummer theory).
 The discussion can be repeated verbatim for the sets
$$H^1_f(\G_p, [\piu]_n)$$
of local Galois cohomology for the group $\G_p:=\Gal(\bQ_p,\Q_p)$.
This local space also admits a map from
$X(\Z_p)$ that fits into a commutative diagram
$$\begin{diagram}
X(\Z_S) & \rTo & X(\Z_p) \\
\dTo & & \dTo \\
H^1_f(\G, [\piu]_n)& \rTo & H^1_f(\G_p, [\piu]_n)
\end{diagram}$$
It comes furthermore with  an analytic description
$$H^1_f(\G_p, [\piu]_n)\simeq [\pi_1^{DR}(X_{\Q_p},b)]_n/F^0$$
provided by {\em $p$-adic Hodge theory } and the {\em De Rham fundamental group
$\pi_1^{DR}(X_{\Q_p},b)$} together with its Hodge filtration $F^{\cdot}$. Thus, eventually, our diagram becomes
$$\begin{diagram}
X(\Z_S) & \rTo & X(\Z_p) \\
\dTo & & \dTo \\
H^1_f(\G, [\piu]_n)& \rTo & [\pi_1^{DR}(X_{\Q_p},b)]_n/F^0
\end{diagram}$$
the effect of which is that
we have replaced the difficult inclusion
$$X(\Z_S) \hra X(\Z_p)$$
with
$$H^1_f(\G, [\piu]_n) \ra [\pi_1^{DR}(X_{\Q_p},b)]_n/F^0,$$
an algebraic map between $\Q_p$-varieties.

It is reasonable  to state a theorem:
\begin{thm}
Let $X$ be a curve and
suppose
$$\dim H^1_f(\G, [\piu]_n) < \dim \pi_1^{DR}(X_{\Q_p},b)_n/F^0$$
for some $n$. Then
$X(\Z_S)$ is finite.
\end{thm}

The proof of the theorem is contained in the following diagram:
$$\begin{diagram}
X(\Z_S) & \rTo & X(\Z_p) \\
\dTo & & \dTo \\
H^1_f(\G, [\piu]_n)& \rTo & [\pi_1^{DR}(X_{\Q_p},b)]_n/F^0\\
& & \dTo^{\a \neq 0}\\
& & \Q_p
\end{diagram}$$
The assumption on dimensions implies that the image of $H^1_f(\G, [\piu]_n)$
inside $\pi_1^{DR}(X_{\Q_p},b)_n/F^0$ is not Zariski dense, and hence, is killed by
some non-zero function $\a$. However, when the function is pulled back to
$X(\Z_p)$ it turns out to be a non-zero linear combination of {\em $p$-adic iterated integrals}
$$\int_b^{x}\b_1 \b_2\cdots \b_m$$
of differential forms $\b_i$ on $X$. This description is the really useful technical input from
$p$-adic Hodge theory.
The point is that such a function can be expanded as a non-vanishing convergent power series on
each $p$-adic disk in $X(\Z_p)$, and hence, has only finitely many zeros. The commutativity of the diagram
is then enough to imply that the function vanishes on $X(\Z_S)$, yielding for us its finiteness.

Some amount of progress has accrued to the program of {\em non-abelian Diophantine geometry}
by way of this theorem, such as new proofs of Diophantine finiteness for hyperbolic curves
of genus zero or one. Furthermore, standard conjectures from the theory
of mixed motives imply that the inequality in the
hypothesis should always hold on hyperbolic curves, insofar one climbs sufficiently
high up on the tower ($n>>0$). One hopes (perhaps in vain) that the milieu of investigation is rich enough
 to  include eventually a broader range of applications,
 such as
a structural understanding of the relationship between Diophantine finiteness and hyperbolically,
and a `non-abelian extension' of the main ideas surrounding the conjecture of
Birch and Swinnerton-Dyer.

In the meanwhile, it is rather interesting to note the key role played by moduli spaces of
principal bundles on $\Spec(\Q)$ such as
$$H^1_f(\G, [\piu]_n).$$
The situation is an appropriate non-abelian complement to the classical use of the Jacobian of a curve,
and the occurrence of related moduli spaces in the Langlands' program. It appears to have been
Andr\'e Weil who first foresaw such possibilities in a remarkable paper of the 1930's, even
with no knowledge of the \'etale topology. This is a point of
considerable historical  interest that will be elaborated upon in a separate
lecture.

\end{document}